%% file: 2002-22.tex
\documentclass{gtart}

\input gtoutput

\lognumber{239}

\volumenumber{6}
\papernumber{22}
\volumeyear{2002}
\pagenumbers{649}{655}
\received{26 Febrary 2002}
\accepted{15 November 2002}
\proposed{David Gabai}
\seconded{Benson Farb, Cameron Gordon}
\published{14 December 2002}

\usepackage{amsmath,amssymb,graphicx}

\let\Bbb\mathbb
\let\bold\mathbf
\theoremstyle{plain}
\newtheorem*{Main}{Main Theorem}
\newtheorem{Thm}{Theorem}

\newtheorem{Lem}[Thm]{Lemma}

\theoremstyle{definition}
\newtheorem*{Def}{Definition}
\newtheorem*{Def1}{Definition of convergence I}
\newtheorem*{Def2}{Definition of convergence II}
\newtheorem{Rem}{Remark}

\newtheorem*{Ex}{Example}

\newcommand{\interior}{^{ \kern-5pt ^\circ}}
\newcommand {\bd}{\partial}

\newcommand {\ol}{\overline}

\newcommand {\N}{{\Bbb N}}

\newcommand {\HH}{{\Bbb H}}

\newcommand {\st}{\text{Stab}}

\newcommand {\cA}{ {\cal A}}
\newcommand {\cB}{ {\cal B}}

\newcommand {\cS}{{\cal S}}

\begin{document}
\title {Convergence groups from subgroups}
\author{Eric L Swenson}
\address{Brigham Young University, Mathematics Department\\292 TMCB, Provo,
UT 84604, USA}

\email{eric@math.byu.edu}

\begin{abstract}
We give  sufficient conditions for a group of homeomorphisms of a
Peano continuum $X$ without cut-points  to be a convergence group.
The condition is that there is a collection of convergence
subgroups whose limit sets ``cut up" $X$ in the correct fashion.
This is closely related to the result in \cite{SWE}.
\end{abstract}

\asciiabstract{ We give sufficient conditions for a group of
homeomorphisms of a Peano continuum X without cut-points to be a
convergence group.  The condition is that there is a collection of
convergence subgroups whose limit sets `cut up' X in the correct
fashion.  This is closely related to the result in [E Swenson, Axial
pairs and convergence groups on S^1, Topology 39 (2000) 229-237].}


\maketitlepage

This paper is part of an approach to the weak
hyperbolization conjecture in the case where the three manifold,
$M$, contains an immersed incompressible surface.  The idea of the
program is to compactify the universal cover of $M$ with a
2--sphere in such a way that $G= \pi_1(M)$ acts on the 2--sphere as
a uniform convergence group, and so by Bowditch, \cite{BOW2}, $G$
is word hyperbolic.

The  application, in this case, is when $X$ is a 2--sphere and the
collection $\cA$ of subgroups are surface subgroups whose limit
sets are circles.

We prove the result in the setting of metric spaces, but it is
also true in the more general setting of Hausdorff spaces, and the
proof goes through with only minor modifications.

\begin{Def1}We say a group $G$ acting by
homeomorphism on a space $X$ acts as a {\em convergence group} on
$X$ if for each sequence of distinct elements of $G$, there is a
subsequence $(g_i)$ and points $p,n \in X$ such that for any
compact $C \not \ni n$ and any neighborhood $U$ of $p$ there is an
$K \in \N$ such that for all $i>K$ , $g_i(C) \subset U$. (See figure
1.)
\end{Def1}

\begin{figure}[ht!]
\centerline{
\includegraphics[width=4in, clip]{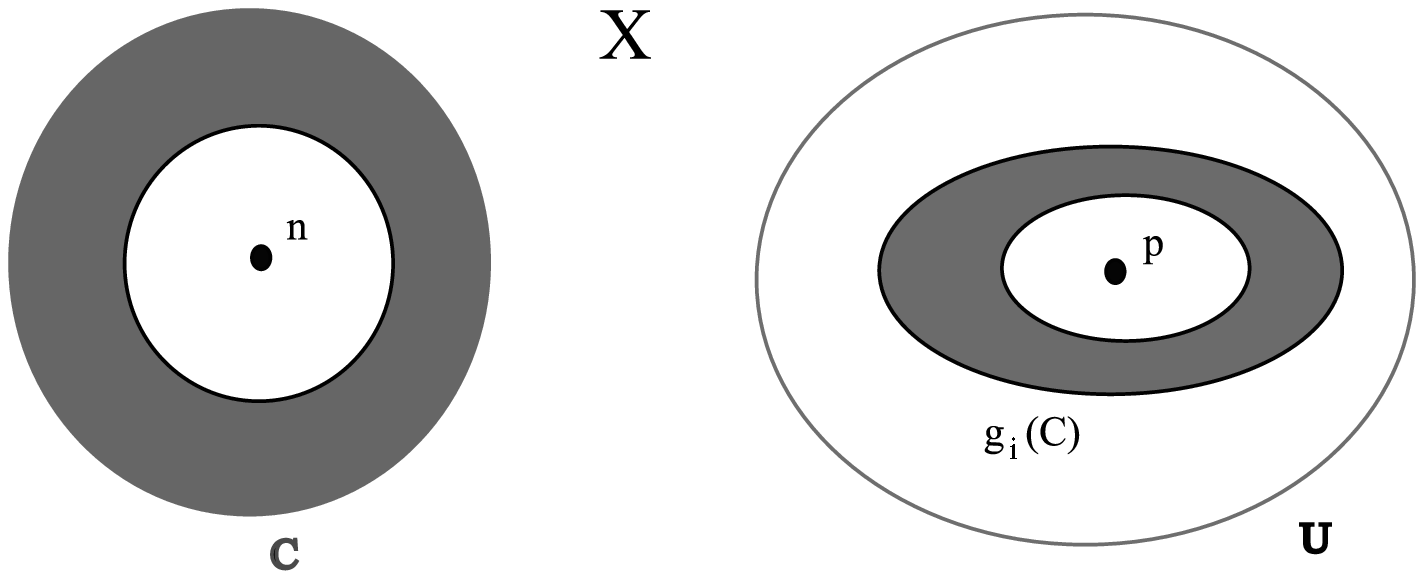}}
\nocolon\caption{}
\end{figure}

\begin{Def} We say that $G$ acts properly on a space
$Y$
if for any compact $C \subset Y$ the set $\{g \in G: g(C) \cap C
\neq \emptyset \}$ is finite
\end{Def}
\begin{Rem}\label{R:pd}
It is easily shown that $G$ fails to act properly on $Y$ if and
only if there exists a sequence $( v_i)$ in $ Y$ and a sequence of
distinct group elements $(g_i)$ in $G$ such that $ v_i \to  v \in
Y$ and $g_i( v_i) \to \hat  v \in Y$
\end{Rem}
\begin{Def2}We say a group $G$ which acts by
homeomorphisms on a space $X$ is a {\em convergence group } on $X$
if the induced action on the space $$\triangle X=\{ T\subset X:\,
|T|=3 \}$$ of distinct triples of $X$ is proper. Notice that
$\triangle X$ is a quotient of a subset of $X^3$ by the action of
$S_3$ (permuting coordinates), and this defines the topology on
$\triangle X$.
\end{Def2}

In \cite{BOW}, these two definitions are shown to be equivalent in
the case where $X$ is compact Hausdorff.
\begin{Def} A {\em Peano} continuum is a compact connected
 locally connect metric space.  In a continuum $Y$, we say
 that $A \subset Y$ {\em separates} the sets $B,C \subset Y$ if $B$
 and $C$ lie in different components of $Y-A$.
 In a Peano continuum $X$, we say that closed sets $A, B \subset X$
  {\em cross} if $A\cap B \neq \emptyset$, $A$ separates points of $B$,
  or $B$ separates points of $A$.
\end{Def}

\begin{Def} Let $X$ be a Peano continuum, and
 $G$ be a group which acts by homeomorphisms on $X$.  If  $\cA$ is a
$G$--invariant collection of closed subsets of $X$, we say that the
pair $(G, \cA)$ is a fine pairing on $X$ if the following
conditions are satisfied:
\begin{enumerate}
\item $\cA$ is { \em cross connected}, that is for any
$A,B \in \cA$ there are  $ A_1, A_2, \dots A_n \in \cA$ with $A_1
=A$, $A_n =B$, where $A_i$ crosses $A_{i+1}$ for $0 < i <n$.
\item $\cA$ is {\em null}.  That is: For any $\epsilon >0$,
the set of elements of $\cA$ with diameter at least $\epsilon$,
$\{A \in \cA: diam(A) > \epsilon \}$ is finite.
\item $\cA$ is {\em fine}.  That is: For any $x,y \in X$
there exists a finite $\cB \subset \cA$ such that $\cup \cB$
separates $x$ from $y$.
\end{enumerate}
\end{Def}
We need the following two results from continuum theory.
\begin{Lem} \label{L:topo} If $M$ is a Peano continuum and a
closed set $A$ does not separate the closed connected set $B$ from
the closed connected set $C$, then there exists a neighborhood $V$
of $A$ which does not separate $B$ from $C$.  (We are of course
assuming that $A \cap (B \cup C) = \emptyset$.)
\end{Lem}
\begin{proof}
For $x,y \in M-A$, define $x \sim y$ if there exist a closed
neighborhood $D$ of $A$ which does not separate $x$ from $y$.
Clearly $B$ is in a single equivalence class as is $C$.  Let  $ b
\in B$. Define $W=\{ x \in M-A: x \sim b \}$. Since $W$ is the
union of components of open sets of $M$, and $M$ is locally
connected, it follows that $W$ is open in $M$. Notice that $\ol W
\subset W \cup A$ and so $W$ is closed in $M-A$. Thus $W$ is a
component of $M-A$ and so $W \supset C$. It follows that there is
a closed neighborhood $V\subset M-(B\cup C)$ of $A$ which does not
separate some point of $C$ from $b$, and so $V$ does not separate
$B$ from $C$.
\end{proof}
\begin{Lem} \label{L:emma} Given $Z$, a Peano continuum, for any $\epsilon >0$
there is $\delta >0$ such that if $A,B$ are crossing closed sets
of $X$ with $diam(A)< \delta$ and $diam(B)< \delta$, then $diam(A
\cup B) < \epsilon$.
\end{Lem}
\begin{proof} Suppose not, then for some $\epsilon >0$
there are sequencea $(A_n)$ and $(B_n)$ of closed subsets of $X$
with $A_n$ crossing $B_n$ for each $n$ , $\,diam(A_n), diam(B_n) <
\frac 1 n\,$, but  $diam(A_n \cup B_n) \ge \epsilon$ for all $n$.
Using compactness, we may assume that $A_n \to a \in X$ and $B_n
\to b \in X$ which implies that $d(a,b) \ge \epsilon$.  Choose
disjoint connected neighborhoods $U $ and $V$ of $a$ and $b$
respectively. For all $n >>0$, $A_n \subset U$ and $B_n \subset
V$, which implies that
 $A_n$ and $B_n$ do not cross.
\end{proof}
\begin{Lem} \label{L:pigion}  Let $(G, \cA)$ be a fine
pairing on $X$ and $\cB \subset \cA$ with $\cB$ finite and cross
connected. If $(g_i)$ is a sequence of elements of $G$ then either
$(g_i(\cup \cB))$ is null, that is $\forall \epsilon >0\, ,
\exists N$ such that $diam(g_i(\cup \cB)) <\epsilon$ whenever $i >
N$, or there exists $B \in \cB$, $A \in \cA$ and a subsequence
$(g_{i_n})$ of $(g_i)$ with $g_{i_n}(B) =A$ for all $n$.
\end{Lem}
\begin{proof} We assume that $(g_i(\cup \cB))$ is not null.
  By Lemma \ref{L:emma} there is $B \in \cB$ with $(g_i(B))$ not null.  Thus
there is $\epsilon >0$ and a subsequence $(g_{i_m})$ with
$diam(g_{i_m}(B)) > \epsilon$  for all $m$.  Since $\cA$ is null,
the set $\{g_{i_m}(B) \}$ is finite and the result follows.
\end{proof}

\begin{Main}  Let $X$ be a Peano continuum without cut
 points, and $(G, \cA)$ is a fine pairing on $X$.  If, for each
$A \in \cA$, $\st (A)$ acts as a convergence group on $X$, then
  $G$ acts as a convergence group on $X$.
\end{Main}
\begin{proof} We will show that $G$ acts properly on
$\triangle X$.  Assume not, then we have a sequence $({\bold
v_i})$ in $\triangle X$ and a sequence of distinct group elements
$(g_i)$ in $G$ such that ${\bold v_i} \to {\bold v} \in \triangle
X$ and $g_i({\bold v_i}) \to  {\bold w} \in \triangle X$.  We may
assume we have  ${\bold v_i} = \{ x_i^1,x_i^2, x_i^3 \}$, ${\bold
v} = \{ x^1, x^2, x^3 \}$ and ${\bold w} =\{ y^1, y^2, y^3 \}$
with $x^j_i \to x^j$ and $g_i(x_i^j) \to y^j$ for $j = 1,2,3$.

For $j=1,2,3$ choose a finite $\cB ^j \subset \cA$ with $x^j \in
V^j$, a component of $X- \cup \cB ^j$.  Using fineness (and
compactness of $X$) we can arrange that $\ol{V^j} \cap \ol{V^k} =
\emptyset$ for $j \neq k$. Now choose a finite sequentially
connected $\cB \subset \cA $ so that $\cB^j \subset \cB$ for $j
=1,2,3$.

We now choose connected neighborhoods $W^j$ of $y^j$ for $j =
1,2,3$ such that $\ol{W^j}\cap \ol{W^k}= \emptyset$ for $j \neq
k$.

We now show that the sequence $(g_i(\cup \cB))$ is null. For if
not, then by Lemma \ref{L:pigion} we may assume that there is a $B
\in \cB$ and $A \in \cA$ such that $g_i(B) = A$ for all $i$.  Thus
$h_i =g_ig^{-1}_1 \in \st (A)$. Letting ${\bold u}_i = g_1({\bold
v}_i)$ we have $h_i({\bold u}_i) \to { \bold w}$ and ${\bold u}_i
\to g_1({\bold v})$.  By remark \ref{R:pd}, this would imply that
$\st(A)$ didn't act properly on  $\triangle X$ which is a
contradiction. Thus the sequence $(g_i(\cup \cB))$ is null and
passing to a subsequence we may assume that $(g_i(\cup \cB)) \to z
\in X$. (That is to say that for any neighborhood $U$ of $z$,
there is an $M$ such that if $i >M$ then $g_i(\cup \cB) \subset
U$.)

Since $\ol{W^1}$, $\ol{W^2}$ and $\ol{W^3}$ are disjoint,
  we may assume $z \not \in \ol{W^2} \cup \ol{W^3}$,
  and so passing to a subsequence we
may assume that $g_i(\cup \cB) \cap W^2 = \emptyset = g_i(\cup
\cB) \cap W^3$. The point $z$ is not a cut point of $X$. By Lemma
\ref{L:topo}, for $i \gg 0$, $g_i(\cup \cB)$ does not separate
$W^2$ from $W^3$. (See figure 2.)

\begin{figure}[ht!]
\centerline{
\includegraphics[width=3in,clip]{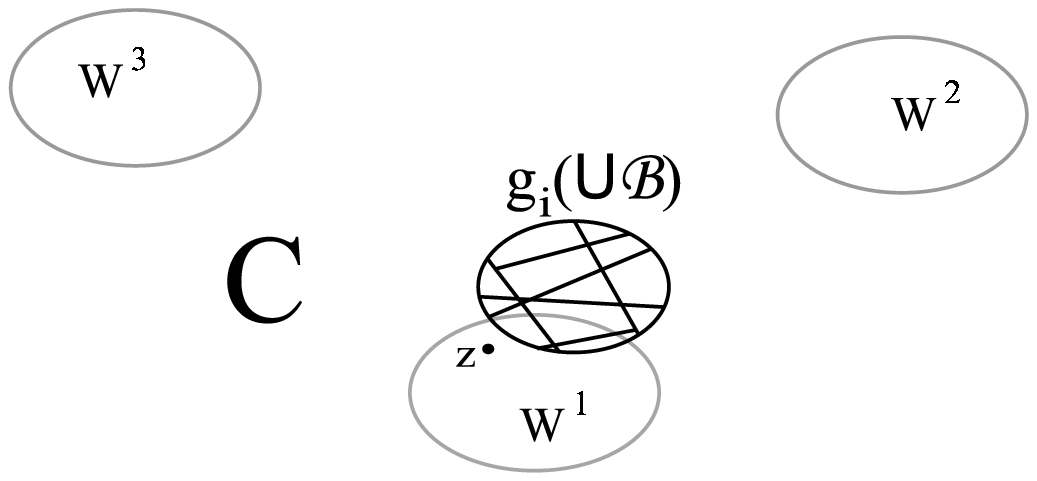}}
\nocolon\caption{}
\end{figure}

Notice that for $i \gg 0$, $g_i(x_i^j) \in g_i(V^j) \cap W^j \neq
\emptyset$.

 Thus for such a (very) sufficiently large $i$, $W^2 \cup W^3$
is contained in a single component $C$ of $X - g_i(\cup \cB)$ and
$g_i(V^2) \cap W^2 \neq \emptyset \neq g_i(V^3) \cap W^3$.  This
implies that $C \subset g_i(V^2)$ and $C \subset g_i(V^3)$
contradicting the fact that $V^2 \cap V^3 = \emptyset$.
\end{proof}
We will show that the theorem is false if the hypothesis that
$\cA$ is fine is removed.
\begin{Ex}
Let $G$ be the Kleinian reflection group generated by reflection
in the sides of a regular right-angle hyperbolic dodecahedron. For
each reflection $r \in G$ let $P_r$ be the plane in $\HH^3$ fixed
by $r$.  By standard results, the centralizer of $r$ acts
cocompactly on $P_r$. Let $H_r =\st(P_r) = \st(\bd P_r)$ where
$\bd P_r$ is the boundary circle of $P_r$ in $S^2 = \bd \HH^3$.
Since $H_r$ acts cocompactly on $P_r$, it follows that the limit
set $\Lambda H_r = \bd P_r$.  Notice that the Main theorem applies
in this setting, as the collection of circles $ \cA = \{\bd P_r
|\, r \in G \text{ a reflection } \}$ is fine.  In fact this
collection of circles gives one of the subdivision rules studied
by Cannon, Floyd and Perry \cite[pages 16--17]{CA-FL-PE}.

We now alter this example so that the collection $\cA$ is no
longer fine, and the action is no longer a convergence action, but
all of the other hypothesis of the Main theorem apply.
 By Baire's theorem, there is a point
$s \in S^2$ such that $s \not \in \bd P_r$ for every reflection $r
\in G$.  Choose a convergence sequence of distinct group elements
$(g_i) \subset G$ with $n \neq s$ such that for any compact $C
\not \ni n$ and any neighborhood $U$ of $s$, $g_i(C) \subset U$
for all $i >>0$.

Now let $\cS$ be the set of translates of $s$ by $G$. For each $t
\in \cS$, we replace the point $t$ by the circle of directions
$\hat t$ at $t$ (so $\hat t$ is a circle of length $2\pi$).
Rigourously let $\hat \cS = \{\hat t | t \in \cS \}$, and  let $X
= (S^2 - \cS)\cup \hat \cS$. We now describe the topology on $X$
and show that $G$ acts on $X$ via homeomorphisms. There is the
obvious "quotient" function $f\co X \to S^2$ given by $f(x) =x$ for
$x \not \in \cup \hat \cS$ and $f(\hat t) = \{t \}$ for $t \in
\cS$.  To define the topology on $X$ we describe a neighborhood
basis of $X$.  For $x \not \in \cup \hat \cS$ and for $U$ a
neighborhood of $f(x)$, $f^{-1}(U)$ will be a neighborhood of $x$.

For $x \in \hat t \in \hat \cS$, let
$\epsilon \in (0, \pi)$. Let $\alpha$ be the geodesic segment
starting at $t$ of length $\epsilon$ in the direction of $x$.
Define
$$U_{\epsilon} = \{y \in S^2: \, 0 < d(t,y) < \epsilon,\text{ and }
\angle_t(\alpha,[t,y]) < \epsilon \}.$$ Now $V_\epsilon =
f^{-1}(U_\epsilon) \cup (x - \epsilon, x+\epsilon)$, where
$(x-\epsilon, x+\epsilon)$ is the open interval of the circle
$\hat t$. (See figure 3.)

\begin{figure}[ht!]
\centerline{
\includegraphics[width=4.5in,clip]{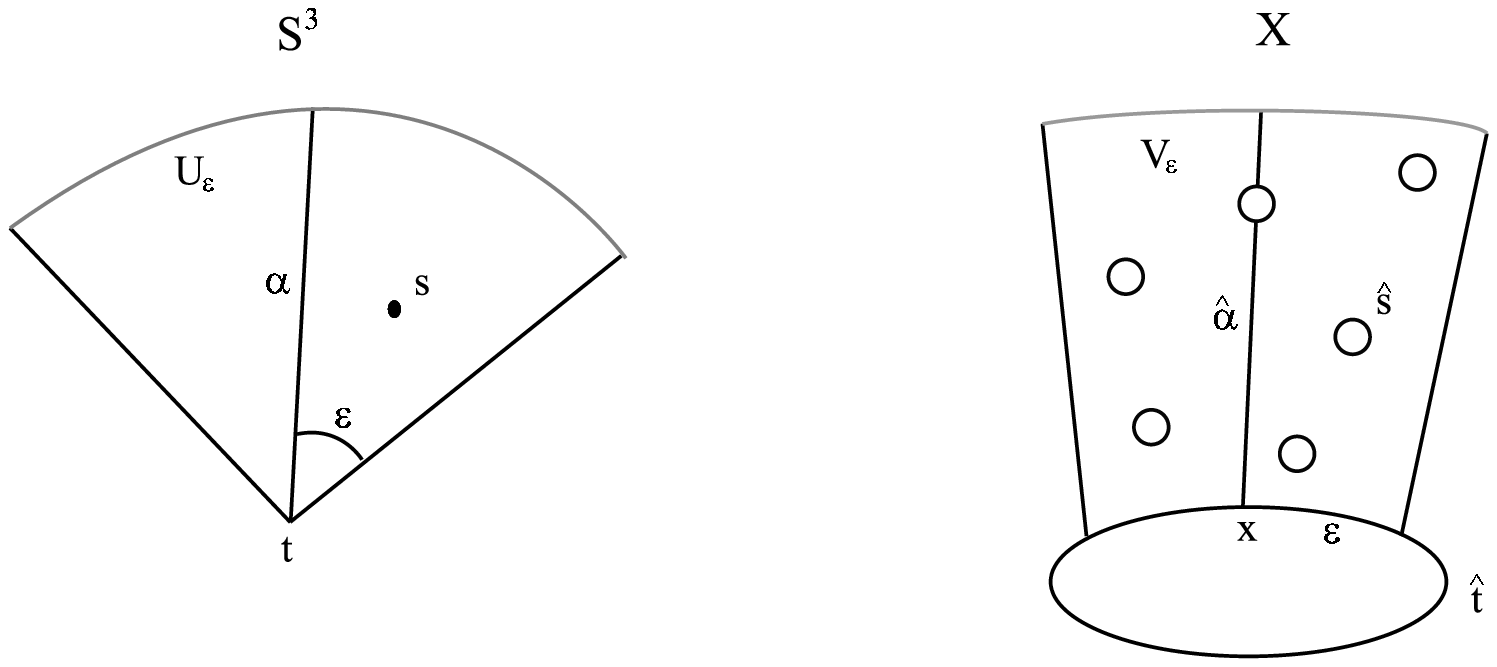}}
\nocolon\caption{}
\end{figure}

This gives a topology on $X$ which is metrizable since $X$ is
regular and has a countable basis ($\cS$ is countable).  Since $G$
acts conformally on $S^2$, $G$ preserves angles, and hence the
action of $G$ on $S^2$ "extends" to an action of $G$ on $X$, and
this action commutes with $f$, ie. for $g \in G$, $g\circ f=
f\circ g$.

 Now for any reflection $r \in G$, by definition of $s$,  $\cS
 \cap \bd P_r = \emptyset$.  Thus $\bd P_r \subset X$ and it
 follows that the
 stabilizer $H_r$ acts as a convergence on $X$.

 However, $G$ does not act as a convergence group on $X$, for our
 original convergence sequence $(g_i)$ will have no convergence
 subsequence.  To see this, let $c \in X$ and
  take a closed  annulus $D \subset S^2$
 separating $n$ from $s$ with $f(c) \not \in D$.
 Notice that by definition of $(g_i)$, $g_i(D) \to s$.
 Let $\hat D = f^{-1}(D)$,  and notice that $\hat D$ is compact.
Clearly for any
 subsequence $(g_{n_i})$  of $(g_i)$, $g_{n_i}(\hat D) \to \hat
 s$.  Since this is converging to an entire circle, it cannot
 converge to a single point, and so $(g_{n_i})$ is not a
 convergence sequence with repelling point $c$.
  Thus $G$ does not act as a convergence
 group on $X$, even though the action of $G$ on $S$ satisfies all of the
 hypothesis of the main theorem except that $\cA$ is not fine on
 $X$.

 We can alter this example slightly by gluing  a disk into each of
 the circles $\hat t$, and obtain a new action of $G$ on $S^2$
 which is not a convergence action (fineness still fails).

\end{Ex}

\end{document}
\bye

%% file: gtoutput.tex
\def\ifplaintex{\expandafter\ifx\csname documentclass\endcsname\relax}

\ifplaintex 
\else
\fi

\expandafter\ifx\csname beginpicture\endcsname\relax
\expandafter\ifx\csname documentclass\endcsname\relax
\input pictex \else
\input prepictex \input pictex \input postpictex \fi\fi

\def\gt{{\mathsurround=0pt\it $\cal G\mskip-2mu$eometry \&\ 
$\cal T\!\!$opology}}        

\def\gtp{{\mathsurround=0pt\it $\cal G\mskip-2mu$eometry \&\ 
$\cal T\!\!$opology $\cal P\!$ublications}}  

\def\lognumber#1{\def\thelognumber{#1}}
\def\volumenumber#1{\def\thevolumenumber{#1}}
\def\papernumber#1{\def\thepapernumber{#1}}
\def\volumeyear#1{\def\thevolumeyear{#1}}

\def\pagenumbers#1#2{\def\startpage{#1}\def\finishpage{#2}}
\def\published#1{\def\publishdate{#1}}
\def\proposed#1{\def\theproposer{#1}}
\def\seconded#1{\def\theseconders{#1}}
\def\received#1{\def\receiveddate{#1}}

\def\accepted#1{\def\accepteddate{#1}}

\long\def\asciiabstract#1{\long\def\theasciiabstract{#1}}

\let\\\par\let\thelognumber\relax
\let\thevolumenumber\relax\let\thepapernumber\relax
\let\thevolumeyear\relax\let\thesamplenumber\relax\let\startpage\relax
\let\finishpage\relax\let\publishdate\relax\let\receiveddate\relax
\let\reviseddate\relax\let\accepteddate\relax\let\theasciititle\relax
\let\theasciiauthors\relax
\let\theasciiabstract\relax
\let\theasciiemail\relax\let\theshortauthors\relax\let\theshorttitle\relax

\long\def\maketitlep{   

\count0=\startpage

\gt\hfill      
\beginpicture
\setcoordinatesystem units <0.33truein, 0.33truein> point at 2.2 0.9
\setplotsymbol ({$\cal G$})
\plotsymbolspacing=9truept
\circulararc 315 degrees from 0 1 center at 0 0
\setplotsymbol ({$\cal T$})
\circulararc 315 degrees from 1 -1 center at 1 0
\endpicture
\break
{\small\ifx\thesamplenumber\relax 
Volume \else Sample
\fi\thevolumenumber\ (\thevolumeyear)
\startpage--\finishpage\nl
Published: \publishdate}
\vglue 0.5truein plus 0.4fil minus 0.1truein

{\parskip=0pt\leftskip 0pt plus 1fil\def\\{\par\smallskip}{\ifplaintex\large
\else\Large\fi\bf\thetitle}\par\medskip}   

\vglue 0pt plus 0.1fil 

{\parskip=0pt\leftskip 0pt plus 1fil\def\\{\par}{\sc\theauthors}
\par\medskip}

\vglue 0pt plus 0.1fil 

{\small\parskip=0pt\let\newline\\
{\leftskip 0pt plus 1fil\def\\{\par}{\sl\theaddress}\par}
\expandafter\ifx\theemail\relax    
\relax\else\vglue 5pt plus 0.02fil minus 2pt\def\\{\stdspace{\rm 
and}\stdspace} 
\cl{Email:\stdspace\tt\theemail}\fi
\ifx\theurl\relax                  
\relax\else\vglue 5pt plus 0.02fil minus 2pt\def\\{\stdspace{\rm 
and}\stdspace}
\cl{URL:\stdspace\tt\theurl}\fi\par}

\vglue 7pt plus 0.3fil minus 3pt

{\bf Abstract}
\vglue 5pt plus 0.1fil minus 2pt

\theabstract

\vglue 7pt plus 0.3fil minus 3pt

{\bf AMS Classification numbers}\quad Primary:\quad \theprimaryclass

Secondary:\quad \thesecondaryclass

\vglue 5pt plus 0.3fil minus 2pt

{\bf Keywords}\quad \thekeywords

\vglue 10pt plus 0.5fil minus 5pt

{\small  Proposed: \theproposer\hfill Received: \receiveddate\nl
Seconded: \theseconders\hfill 
\ifx\reviseddate\relax                         
Accepted: \accepteddate                        
\else
Revised: \reviseddate                          
\fi}
\eject
}       

\let\maketitlepage\maketitlep
\let\maketitle\maketitlepage

\font\phead=cmsl9 scaled 950
\font=cmsl9 scaled 1050
\font\pnum=cmbx10 scaled 913
\font\lnum=cmbx10 
\font\pfoot=cmsl9 scaled 950
\font=cmsl9 scaled 1050
\ifplaintex
\headline{\vbox to 0pt{\vskip -4.5mm\line{\small\phead\ifnum
\count0=\startpage ISSN 1364-0380 (on line)
1465-3060 (printed) \hfill {\pnum\folio}\else\ifodd\count0\def\\{ }%
\ifx\theshorttitle\relax\thetitle\else\theshorttitle\fi\hfill{\pnum\folio}
\else\def\\{ and }{\pnum\folio}\hfill\ifx\theshortauthors\relax\theauthors
\else\theshortauthors\fi\fi\fi}\vss}}
\footline{\vbox to 0pt{\vglue 0mm\line{\small\pfoot\ifnum\count0=\startpage
\copyright\ \gtp\hfill\else
\gt, Volume \thevolumenumber\ (\thevolumeyear)\hfill\fi}\vss
}}
\else
\makeatletter
\def\@oddhead{{\small\lhead\ifnum\count0=\startpage ISSN 1364-0380 (on line)
1465-3060 (printed) \hfill {\lnum\number\count0}\else\ifodd\count0
\def\\{ }\ifx\theshorttitle\relax \thetitle \else\theshorttitle\fi\hfill
{\lnum\number\count0}\else\def\\{ and }{\lnum\number\count0}
\hfill\ifx\theshortauthors\relax 
\theauthors\else\theshortauthors\fi\fi\fi}}\def\@evenhead{\@oddhead}
\def\@oddfoot{\small\lfoot\ifnum\count0=\startpage\copyright\ \gtp\hfill\else
\gt, Volume \thevolumenumber\ (\thevolumeyear)\hfill\fi}
\def\@evenfoot{\@oddfoot}
\makeatother
\fi

\newwrite\gtoutfile
\long\gdef\makeheadfile{  
{\def\\{, }\def\s{ }
\immediate\openout\gtoutfile head.xxx
\immediate\write\gtoutfile{To: math@arxiv.org}
\immediate\write\gtoutfile{Subject: put or rep NNNNN:pppp}
\immediate\write\gtoutfile{--text follows this line--}
\immediate\write\gtoutfile{Proxy-for: \ifx\theasciiauthors\relax
\theauthors\else\theasciiauthors\fi\s<\ifx\theasciiemail\relax\theemail\else\theasciiemail\fi>}
\immediate\write\gtoutfile{\noexpand\\}
\immediate\write\gtoutfile{Authors: \ifx\theasciiauthors\relax
\theauthors\else\theasciiauthors\fi}
{\def\\{ }\immediate\write\gtoutfile{Title: \ifx\theasciititle\relax
\thetitle\else\theasciititle\fi}}
\immediate\write\gtoutfile{Subj-class: GT or SG or MG etc}
\immediate\write\gtoutfile{MSC-class: \theprimaryclass\ifx\thesecondaryclass\relax\else, \thesecondaryclass\fi}
\immediate\write\gtoutfile{Journal-ref: Geom. Topol. \thevolumenumber
(\thevolumeyear) \startpage-\finishpage}
\immediate\write\gtoutfile{Comments: Published by Geometry and Topology at}
\immediate\write\gtoutfile{\s\s http://www.maths.warwick.ac.uk/gt/GTVol\thevolumenumber/paper\thepapernumber.abs.html}
\immediate\write\gtoutfile{\noexpand\\}
\immediate\write\gtoutfile{}
\ifx\theasciiabstract\relax
\immediate\write\gtoutfile{\theabstract}\else
\immediate\write\gtoutfile{\theasciiabstract}\fi
\immediate\write\gtoutfile{}
\immediate\write\gtoutfile{\noexpand\\}
\immediate\write\gtoutfile{}
\immediate\closeout\gtoutfile}}  

\def\maketitlepage{\maketitlep\makeheadfile}
\let\maketitle\maketitlepage